\documentclass[12pt]{article}
\usepackage{lineno,hyperref}
\usepackage{array}
\usepackage[all]{xy}
\usepackage{xcolor}
\usepackage{amsxtra, microtype}
\usepackage{amssymb,amsmath,amsthm,stmaryrd,latexsym,wasysym}
\usepackage{amstext}
\usepackage{dsfont}
\usepackage{graphicx}
\usepackage{leftidx}
\usepackage[version=3]{mhchem}

\newtheorem{theorem}{Theorem}[section]

\newtheorem{definition}[theorem]{Definition}

\newtheorem{lemma}[theorem]{Lemma}

\newtheorem{proposition}[theorem]{Proposition}
\newtheorem{remark}[theorem]{Remark}

\numberwithin{equation}{section}
\newcommand{\duer}{\mathbin{\raisebox{3pt}{\varhexstar}\kern-3.70pt{\rule{0.15pt}{4pt}}}\,}
\usepackage{cite}

\usepackage{enumitem}  
\usepackage{calc}
\setlist{labelindent=1pt,itemsep=.5em}
\setlist[itemize]{leftmargin=1.2cm}
\setlist[enumerate]{itemindent=0em,leftmargin=1.2cm}
\setlist[enumerate,1]{label={\upshape(\roman*)}}

\allowdisplaybreaks

\makeatletter
\newcommand{\subjclass}[2][2010]{%
  \let\@oldtitle\@title%
  \gdef\@title{\@oldtitle\footnotetext{#1 \emph{Mathematics subject classification}: #2}}%
}
\newcommand{\keywords}[1]{%
  \let\@@oldtitle\@title%
  \gdef\@title{\@@oldtitle\footnotetext{\emph{Keywords}: #1.}}%
}
\makeatother

\title{Simply complete hom-Lie superalgebras and decomposition of complete hom-Lie superalgebras}%

\author{M. R. Farhangdoost$^1$\footnote{Corresponding author}, A. R. Attari Polsangi$^1$, S. Silvestrov$^2$ \\
\footnotesize $^1$Department of Mathematics, College of Sciences, Shiraz University, \\
\footnotesize P.O. Box 71457- 44776, Shiraz, Iran \\
\footnotesize \text{farhang@shirazu.ac.ir}   \\
\footnotesize $^2$Division of Mathematics and Physics, School of Education, Culture and Communication, \\
\footnotesize M\"{a}lardalen University, Box 883, 72123 V\"{a}ster{\aa}s, Sweden \\
\footnotesize \text{sergei.silvestrov@mdh.se}}

\subjclass[2020]{17B61, 17D30 (Primary) 17B65, 17B68, 17B70 (Secondary)}
\keywords{hom-Lie superalgebra, complete hom-Lie superalgebra, simple hom-Lie superalgebra algebras}

\date{\today}

\begin{document}

\maketitle

\begin{abstract}
Complete hom-Lie superalgebra are considered and some equivalent conditions for a hom-Lie superalgebra to be a complete hom-Lie superalgebra are established. In particular, the relation between decomposition and completeness for a hom-Lie superalgebra is described. Moreover, some conditions that the set of  $\alpha^{s}$-derivations of a hom-Lie superalgebra to be complete and simply complete are obtained.
\end{abstract}



\section{Introduction}
Hom-Lie algebras and quasi-hom-Lie algebras were introduced first by Hartwig, Larsson, and Silvestrov in 2003 in \cite{HLSPrepr2003JA2006:deformLiealgsigmaderiv} devoted to a general method for construction of deformations and discretizations of Lie algebras of vector fields and deformations of Witt and Virasoro type algebras based on general twisted derivations ($\sigma$-derivations) obeying twisted Leibniz rule, and motivated also by the examples of $q$-deformed Jacobi identities in $q$-deformations of Witt and Visaroro algebras and in related $q$-deformed algebras discovered in 1990'th in string theory, vertex models of conformal field theory, quantum field theory and quantum mechanics, and $q$-deformed differential calculi and $q$-deformed homological algebra
\cite{aizawasato199091:qdefViralgcenext,ChaiElinPop1990:qconfalgcentext,ChaiIsLukPopPresn1991:Viralgconfdim,ChaiKuLuk,ChaiPopPres,
daskaloyannis1992generalized,Hu1992:qWittalgqLie,Kassel92,LiuKQuantumCentExt,LiuKQ1992:CharQuantWittAlg,LiuKQPhDthesis}.
In 2005, Larsson and Silvestrov introduced quasi-Lie and quasi-Leibniz algebras in \cite{LarssonSilvestrov2005:QuasiLiealgebras} and graded color quasi-Lie and graded color quasi-Leibniz algebras in \cite{LarssonSilvestrov2005:GradedquasiLiealgebras}
incorporating within the same framework thehom-Lie algebras and quasi-hom-Lie algebras, the color hom-Lie algebras and hom-Lie superalgebras, quasi-hom-Lie color algebras, quasi-hom-Lie superalgebras, quasi-Leibniz algebras and graded color quasi-Leibniz algebras.
The central extensions and cocycle conditions have been first considered for quasi-hom-Lie algebras and hom-Lie algebras in
\cite{HLSPrepr2003JA2006:deformLiealgsigmaderiv,LarssonSilvJA2005:QuasiHomLieCentExt2cocyid} and for graded color quasi-hom-Lie algebras in \cite{SigurdssonSilvestrov2009:colorHomLiealgebrascentralext}.

In quasi-Lie algebras, the skew-symmetry and the Jacobi identity are twisted by deforming twisting linear maps, with the Jacobi identity in quasi-Lie and quasi-hom-Lie algebras in general containing six twisted triple bracket terms. In hom-Lie algebras, the bilinear product satisfies the non-twisted skew-symmetry property as in Lie algebras, and the hom-Lie algebras Jacobi identity has three terms twisted by a single linear map. Lie algebras are a special case of hom-Lie algebras when the twisting linear map is the identity map. For other twisting linear maps however the hom-Lie algebras are different and in many ways richer algebraic structures with classifications, deformations, representations, morphisms, derivations and homological structures in the fundamental ways dependent on joint properties of the twisting map and bilinear product which are in the intricate way interlinked by hom-Jacobi identity. Hom-Lie admissible algebras have been considered first in \cite{MakhoufSilvestrov:Prep2006JGLTA2008:homstructure}, where
the hom-associative algebras and more general $G$-hom-associative algebras including the Hom-Vinberg algebras (hom-left symmetric algebras), hom-pre-Lie algebras (hom-right symmetric algebras), and some other new Hom-algebra structures have been introduced and shown to be Hom-Lie admissible, in the sense that the operation of commutator as new product in these hom-algebras structures yields hom-Lie algebras. Furthermore, in \cite{MakhoufSilvestrov:Prep2006JGLTA2008:homstructure}, flexible hom-algebras have been introduced and connections to hom-algebra generalizations of derivations and of adjoint derivations maps have been considered, investigations of the classification problems for hom-Lie algebras have been initiated with constriction of families of the low-dimensional hom-Lie algebras.

The Hom-Lie superalgebras and the more general color quasi-Lie algebras provide new general parametric families of non-associative structures, extending and interpolating on the fundamental level of defining identities between the Lie algebras, Lie superalgebras, color Lie algebras and some other important related non-associative structures, their deformations and discritizations,
in the special interesting ways which may be useful for unification of models of classical and quantum physics, geometry and symmetry analysis, and also in algebraic analysis of computational methods and algorithms involving linear and non-linear descretizations of differential and integral calculi.
Investigation of color hom-Lie algebras and hom-Lie superalgebras and $n$-ary generalizations have been further expanded recently in
\cite{AbAmMakh:HomaltHomMalcHomJordSuperalg,
AbdaouiAmmarMakhloufCohhomLiecolalg2015,
AbramovSilvestrov:3homLiealgsigmaderivINvol,
AmmarMakhloufJA2010:homliesuperaladmsuperalg,
AmmarMakhloufSaadaoui2013:CohlgHomLiesupqdefWittSup,
AmmarMakhloufSilv:TernaryqVirasoroHomNambuLie,
AmAyMabMakh:QuadrColHomLieAlgs,
ArmakanFarhangdoost:2017IJGMMP:GeomaspectsextHomLiesuperalgs,
ArmakanFarhSilv2021:ndKilformsHomLiesuper,
armakanrazavicomalg2020:completehomliesuper,
ArmakanSilv2020:envelalgcertaintypescolorHomLie,
ArmakanSilvFarh:envelopalgcolhomLiealg,
ArmakanSilvFarh20172019:exthomLiecoloralg,
ArmakanSilv:colorHomLieHomLiebnOmniHomLie,
akms:ternary,
ams:ternary,
ArnlindMakhloufSilvnaryHomLieNambuJMP2011,
AtMaSi:GenNambuAlg,
BeitesKaygorodovPopov2018:GenDermultnaryHomOmegacolalg,
Bakayoko2014:ModulescolorHomPoisson,
BakayokoDialo2015:genHomalgebrastr,
BakyokoSilvestrov:Homleftsymmetriccolordialgebras,
BakyokoSilvestrov:MultiplicnHomLiecoloralg,
BakayokoToure2019:genHomalgebrastr,
CaoChen2012:SplitregularhomLiecoloralg,
GuanChenSun:HomLieSuperalgebras,
KitouniMakhloufSilvestrov:nhominduced,
kms:narygenBiHomLieBiHomassalgebras2020,
kms:solvnilpnhomlie2020,
MabroukNcibSilvestrov2020:GenDerRotaBaxterOpsnaryHomNambuSuperalgs,
Makhlouf2010:paradigmnonasshomalghomsuper,
MandalMishra:HomGerstenhaberHomLiealgebroids,
MishraSilvestrov:SpringerAAS2020HomGerstenhalgsHomLiealgds,
SigSilv:CzechJP2006:GradedquasiLiealgWitt,
SigurdssonSilvestrov2009:colorHomLiealgebrascentralext,
SilvestrovParadigmQLieQhomLie2007,
WangZhangWei2015:HomLeibnizsuperalg,
Yuan2012:HomLiecoloralgstr,
ZhouChenMa:GenDerHomLiesuper}.

In \cite{chun1996on,goto1975faithful}, the complete Lie superalgebras were introduced and studied.
Recently the notion of compact hom-Lie superalgebra was introduced in \cite{armakanrazavicomalg2020:completehomliesuper}.
In this article, complete hom-Lie superalgebras are considered and equivalent conditions for a hom-Lie superalgebra to be a complete hom-Lie superalgebra are established. In particular, the relation between decomposition and completness for a hom-Lie superalgebra is described. Moreover, some conditions for the set of  $\alpha^{s}$-derivations of a hom-Lie superalgebra to be complete and simply complete are obtained. In Section \ref{sec:prelimhomliesuperalg}, some necessary notations and useful definitions and properties of hom-Lie superalgebras are reviewed. In Section \ref{sec:complhomliesuperalg}, the notion of a complete hom-Lie superalgebra is presented, and the equivalent conditions for the completeness of $\mathfrak{g}_0$ and $\mathfrak{g}$ are studied. Then conditions for a hom-Lie superalgebra to be complete are considered by using the notion of holomorph hom-Lie superalgebras and hom-ideals. After that simply complete hom-Lie superalgebras are defined and equivalence of a hom-Lie superalgebra being simply complete or  indecomposable is investigated. Finally we discuss the conditions for the $Der_{\alpha^{s+1}}(\mathfrak{g})$ to be complete and simply complete.

\section{Preliminaries on hom-Lie superalgebras and their representation and derivations}
\label{sec:prelimhomliesuperalg}
Throughout this article, all linear spaces are assumed to be over a field $\mathbb{K}$ of characteristic different from $2$. A linear space $V$ is said to be a $G$-graded by an abelian group $G$ if, there exists a family $\{V_g\}_{g\in G}$ of linear subspaces of $V$ such that $V=\bigoplus\limits_{g\in G} V_g.$
The elements of $V_g$ are said to be homogeneous of degree $g\in G$.
The set of all homogeneous elements of $V$ is denoted $\mathcal{H}(V)= \bigcup\limits_{g\in G} V_g$.
A linear mapping $f : V\rightarrow V'$ of two $G$-graded linear spaces $V=\bigoplus\limits_{g\in G} V_g$ and $V'=\bigoplus\limits_{g\in G} V'_g$ is called homogeneous of degree $d$ if $f(V_g)\subseteq  V'_{g+d}$, for all $g\in G.$
Homogeneous linear maps of degree zero, $f(V_g)\subseteq V'_{g}$ for any $g\in G$, are also called even.
An algebra $(A, \cdot)$ is said to be $G$-graded if its underlying linear space is $G$-graded, $A=\bigoplus\limits_{g\in G}A_g$, and moreover $A_g\cdot A_h\subseteq A_{g+h}$, for all $g, h\in G.$
A homomorphism $f : A\rightarrow A'$ of $G$-graded algebras $A$ and $A'$
is an algebra morphism which is even (degree $0_G$).
In $\mathbb{Z}_2$-graded linear spaces $A=A_0\oplus A_1$, the elements of $A_j$ are homogeneous of degree (parity) $j\in\mathbb{Z}_2$, and the set of all homogeneous elements is $\mathcal{H}(A)=A_0\cup A_1$. The parity of a homogeneous element $x\in \mathcal{H}(A)$ is denoted $|x|$.

Hom-superalgebras are triples $(A, \mu, \alpha)$ in which $A=A_0\oplus A_1$ is a $\mathbb{Z}_2$-graded linear space ($\mathbb{K}$-superspace), $\mu : A\times A\rightarrow A$ is an even bilinear map, and $\alpha : A\rightarrow A$ is an even linear map. An even linear map $f : A \rightarrow A'$ is said to be a weak morphism of hom-superalgebras if
$f\circ\mu=\mu\circ(f\otimes f),$ and a morphism of hom-superalgebras if moreover $f\circ\alpha=\alpha'\circ f$.

In any hom-superalgebra $(A=A_0\oplus A_1, \mu, \alpha)$,
\begin{align*}
\mu(A_0, A_0)\subseteq A_0, \quad \mu(A_1, A_0)\subseteq A_1, \\
\mu(A_0, A_1)\subseteq A_1, \quad \mu(A_1, A_1)\subseteq A_0.
\end{align*}
Hom-subalgebras (graded hom-subalgebras) of hom-superalgebra $(A,\mu,\alpha)$ are defined as  $\mathds{Z}_{2}$-graded linear subspaces $I =  (I\cap A_{0}) \oplus (I\cap A_{1})$ of $A$ closed  under both $\alpha$ and $\mu$, that is  $\alpha(I) \subseteq I$ and $\mu(I,I) \subseteq I$.

Hom-subalgebra $I$ is called a left hom-ideal of the hom-superalgebra $A$, if $\mu(A, I) \subseteq I$, and the right hom-ideal of $A$ if $\mu(I,A) \subseteq I$. If $I$ in $A$ is both left and right hom-ideal, then it is called hom-ideal or two-sided hom-ideal, and notation $I \triangleleft A$ is used to indicate that. Hom-ideals $I \triangleleft A$ of a hom-superalgebra $(A=A_0\oplus A_1, \mu, \alpha)$, as  $\mathds{Z}_{2}$-graded hom-subalgebras with homogeneous components $I_0 =  (I\cap A_{0})$
and  $(I\cap A_{1})$, satisfy
\begin{align*}
\mu(I_0, A_0)\cup \mu(A_0,I_0) \in I_0, \\
\mu(I_0, A_1)\cup \mu(A_1,I_0) \in I_1 \\
\mu(I_1, A_0)\cup \mu(A_0,I_1) \in I_1  \\
\mu(I_1, A_1)\cup \mu(A_1,I_1) \in I_0.
\end{align*}
If $I \triangleleft A$ and moreover $\alpha(I_0)\subseteq I_0$, then
$(I_0,\mu,\alpha_0)$ is a hom-subalgebra of the hom-algebra $(A_0,\mu,\alpha_0)$, where $\alpha_j:A_j\rightarrow A_j$ are the restrictions of the even linear map  $\alpha:A\rightarrow A$ to homogeneous subspaces $A_j$ for $j\in \mathds{Z}_{2}$. However,  $(I_1,\mu,\alpha_1)$ is a hom-Lie subalgebra of
the hom-algebra $(A_1,\mu,\alpha_1)$ if and only if $\alpha(I_1)\subseteq I_1$ and $\mu(I_1,I_1)=\{0\}$, since $\mu(I_1,I_1)\subseteq A_0 \cap I_1 \subseteq A_0\cap A_1 = \{0\}$.

In any hom-superalgebra $(A=A_0\oplus A_1, \mu, \alpha)$,
hom-associator of $A$ is the even trilinear map given by  $as_{\alpha,\mu}=\mu\circ(\mu\otimes\alpha-\alpha\otimes\mu):A\times A\times A\rightarrow A,$
acting on elements as
$as_{\alpha,\mu}(x, y, z)=\mu(\mu(x, y), \alpha(z))-\mu(\alpha(x), \mu(y, z)),$
or $as_{\alpha,\mu}(x, y, z)=(xy)\alpha(z)-\alpha(x)(yz)$ in juxtaposition notation $xy=\mu(x, y)$.
Since $|as_{\alpha,\mu}(x,y,z))| =  |x|+|y|+|z|$ for $x, y, z \in \mathcal{H}(A)=A_0\cup A_1$ in any hom-superalgebra $(A=A_0\oplus A_1,\mu, \alpha),$
\begin{align}
as_{\alpha,\mu}(A_0, A_0, A_0)\subseteq A_0, & \quad
as_{\alpha,\mu}(A_1, A_1, A_0)\subseteq A_0, \\
as_{\alpha,\mu}(A_1, A_0, A_1)\subseteq A_0, & \quad
as_{\alpha,\mu}(A_0, A_1, A_1)\subseteq A_0, \\
as_{\alpha,\mu}(A_1, A_0, A_0)\subseteq A_1, & \quad
as_{\alpha,\mu}(A_0, A_1, A_0)\subseteq A_1, \\
as_{\alpha,\mu}(A_0, A_0, A_1)\subseteq A_1, & \quad
as_{\alpha,\mu}(A_1, A_1, A_1)\subseteq A_1.
\end{align}
Hom-ideals $I \triangleleft A$ of a hom-superalgebra $(A=A_0\oplus A_1, \mu, \alpha)$, as  $\mathds{Z}_{2}$-graded hom-subalgebras with homogeneous components $I_0 =  (I\cap A_{0})$
and  $I_1=(I\cap A_{1})$, satisfy
\begin{align*}
& as_{\alpha,\mu}(A_j, I_0, A_j)\subseteq I_0, \quad j\in \mathds{Z}_{2} \\
& as_{\alpha,\mu}(A_0, I_1, A_1)\cup as_{\alpha,\mu}(A_1, I_1, A_0) \subseteq I_0, \\
& as_{\alpha,\mu}(A_1, I_1, A_1) \subseteq I_1.
\end{align*}
In particular, for each $j= \{0,1\}\in \mathds{Z}_{2}$, if $\alpha(I_j) \subseteq I_j$, then $I_j$ is closed under ternary trilinear product defined by the hom-asociator $as_{\alpha,\mu}$ which together with $\alpha_j$ define then the structure of ternary hom-algebra on $I_j$. In particular, when
$\alpha(I_0) \subseteq I_0$ and $\alpha(I_1) \subseteq I_1$, both $(I_0,as_{\alpha,\mu},\alpha_0)$ and $(I_1,as_{\alpha,\mu},\alpha_1)$ become ternary hom-algebras at the same time.

\begin{definition}[\cite{HLSPrepr2003JA2006:deformLiealgsigmaderiv,MakhoufSilvestrov:Prep2006JGLTA2008:homstructure}]
\label{def:homlie}	
Hom-Lie algebras are triples $(\mathfrak{g},[.,.],\alpha)$, where $\mathfrak{g}$ is a linear space,
	$[.,.] : \mathfrak{g} \times  \mathfrak{g} \to$  $\mathfrak{g}$ is a bilinear map
	and
	$\alpha: \mathfrak{g} \to \mathfrak{g}$ is a linear map
	satisfying for all $x, y, z \in \mathfrak{g}$,
	\begin{align}
[x, y]&=-[y,x] & \text{Skew-symmetry} \label{skewsymmetry} \\
[\alpha (x), [y, z]]+ [\alpha (y), [z, x]] &+ [\alpha (z), [x, y]] = 0, & \text{Hom-Lie Jacobi identity}
\label{homliejacobi}
	\end{align}
	\begin{enumerate}[label=\upshape{(\roman*)},leftmargin=30pt]
\item Hom-Lie algebra is called a multiplicative hom-Lie algebra if $\alpha$ is an algebra morphism, $\alpha([.,.]) = ([\alpha(.),\alpha(.)])$, meaning that $\alpha([x,y])=[\alpha(x),\alpha(y)]$ for any $x, y \in \mathfrak{g}$.
\item 	Multiplicative hom-Lie algebra is called regular, if $\alpha$ is an automorphism.
	\end{enumerate}
\end{definition}

From the point of view of Hom-algebras, Lie algebras are a special subclass of Hom-Lie algebras obtained when  $\alpha = id$ in Definition \ref{def:homlie}.

\begin{definition}[\cite{AmmarMakhloufJA2010:homliesuperaladmsuperalg,LarssonSilvestrov2005:GradedquasiLiealgebras}] \label{def:homliesuper} Hom-Lie superalgebras are tripples $(\mathfrak{g},[.,.],\alpha)$ which consist of $\mathds{Z}_{2}$-graded linear space $\mathfrak{g} = \mathfrak{g}_{0} \oplus \mathfrak{g}_{1}$, an even bilinear map  $[.,.] : \mathfrak{g} \times \mathfrak{g} \to \mathfrak{g}$ and an even linear map $\alpha : \mathfrak{g} \to \mathfrak{g}$ satisfying the super skew-symmetry and hom-Lie super Jacobi identities for homogeneous elements $x, y, z \in \mathcal{H}(\mathfrak{g})$,
\begin{align}
[x,y] =-(-1)^{|x||y|}[y,x], \quad \quad \quad  \text{Super skew-symmetry}   \label{superskewsymmetry} \\
(-1)^{|x||z|}[\alpha(x),[y,z]] +(-1)^{|y||x|}[\alpha(y),[z,x]] +(-1)^{|z||y|}[\alpha(z),[x,y]]=0. \label{homliesuperjacobi}  \\
\text{Super Hom-Jacobi identity} \nonumber
\end{align}
	\begin{enumerate}[label=\upshape{(\roman*)},leftmargin=30pt]
		\item Hom-Lie superalgebra is called multiplicative Hom-Lie superalgebra, if $\alpha$ is an algebra morphism, $\alpha([x,y])=[\alpha(x),\alpha(y)]$ for any $x, y \in \mathfrak{g}$. 	
		\item Multiplicative hom-Lie superalgebra is called regular, if $\alpha$ is an algebra automorphism.
	\end{enumerate}
\end{definition}

In skew-symmetric hom-superalgebras, the super hom-Jacobi identity can be presented equivalently in the form of super hom-Leibniz rule for the maps $ad_x=[x,.]:\mathfrak{g} \rightarrow \mathfrak{g},$
\begin{equation}
[\alpha(x),[y,z]]=[[x,y],\alpha(z)]+(-1)^{|x||y|}[\alpha(y),[x,z]] \label{homliesuperjacobi:leibniz}
\end{equation}
since, by super skew-symmetry \eqref{superskewsymmetry}, the following equalities are equivalent:
\begin{align*}
& [\alpha(x),[y,z]] = [[x,y],\alpha(z)]+(-1)^{|x||y|}[\alpha(y),[x,z]], \\
& [\alpha(x),[y,z]]-[[x,y],\alpha(z)]-(-1)^{|x||y|}[\alpha(y),[x,z]]=0, \\
& [\alpha(x),[y,z]]+(-1)^{|z|(|x|+|y|)}[\alpha(z),[x,y]]-(-1)^{|x||y|}[\alpha(y),[x,z]]=0, \\
& [\alpha(x),[y,z]]+(-1)^{|z|(|x|+|y|)}[\alpha(z),[x,y]]-(-1)^{|x||y|}[\alpha(y),[x,z]]=0, \\
& [\alpha(x),[y,z]]+(-1)^{|z|(|x|+|y|)}[\alpha(z),[x,y]]+(-1)^{|x||y|}(-1)^{|z||x|}[\alpha(y),[z,x]]=0, \\
& (-1)^{|z||x|}[\alpha(x),[y,z]]+ (-1)^{|z||x|}(-1)^{|z|(|x|+|y|)}[\alpha(z),[x,y]]+(-1)^{|x||y|}[\alpha(y),[z,x]]=0, \\
& (-1)^{|x||z|}[\alpha(x),[y,z]]+(-1)^{|z||y|}[\alpha(z),[x,y]]+(-1)^{|y||x|}[\alpha(y),[z,x]]=0, \\
& (-1)^{|x||z|}[\alpha(x),[y,z]]+(-1)^{|z||y|}[\alpha(z),[x,y]]+(-1)^{|y||x|}[\alpha(y),[z,x]]=0.
\end{align*}

\begin{remark}
If skew-symmetry \eqref{skewsymmetry} does not hold, then \eqref{homliesuperjacobi} and
\eqref{homliesuperjacobi:leibniz} are not necessarily equivalent, defining different Hom-superalgebra structures.
The Hom-superalgebras defined by just super algebras identity \eqref{homliesuperjacobi:leibniz} without requiring super hom-skew-symmetry on homogeneous elements are Leibniz Hom-superalgebras, a special class of general $\Gamma$-graded quasi-Leibniz algebras (color quasi-Leibniz algebras) first introduced in \textrm{\cite{LarssonSilvestrov2005:QuasiLiealgebras, LarssonSilvestrov2005:GradedquasiLiealgebras}}.
\end{remark}
\begin{remark}
In any hom-Lie superalgebra, $(\mathfrak{g}_{0},[.,.],\alpha)$ is a hom-Lie algebra since
$[\mathfrak{g}_{0},\mathfrak{g}_{0}]\in \mathfrak{g}_{0}$ and $\alpha(\mathfrak{g}_{0})\in \mathfrak{g}_{0}$ and $(-1)^{|a||b|}= (-1)^{0}=1$ for $a,b \in \mathfrak{g}_{0}$. Thus,
hom-Lie algebras can be also seen as special class of hom-Lie superalgebras when $\mathfrak{g}_{1}=\{0\}$.
\end{remark}
\begin{remark}
From the point of view of Hom-superalgebras, Lie superalgebras is an important subclass of Hom-Lie superalgebras obtained when $\alpha = id$ in Definition \ref{def:homliesuper}.
Namely, when $\alpha=Id$, Definition \ref{def:homliesuper} becomes the definition of Lie superalgebras \textrm{\cite{Berezin83:algananticomvar,Berezin87:introsuperanalysis,kac1977lie,kac1977lie}} as  $\mathds{Z}_{2}$-graded linear spaces $\mathfrak{g} = \mathfrak{g}_{0} \oplus \mathfrak{g}_{1}$, with a graded Lie bracket $[.,.] : \mathfrak{g} \times \mathfrak{g} \to \mathfrak{g}$ of degree zero, that is $[.,.]$ is a bilinear map obeying
$[\mathfrak{g}_{i}, \mathfrak{g}_{j}] \subset \mathfrak{g}_{{i+j}(mod2)},$
and for homogeneous $x, y, z \in \mathcal{H}(\mathfrak{g})$,
\begin{align}
[x,y]=-(-1)^{|x||y|}[y,x],  \quad \quad \quad \text{Super skew-symmetry}   \\
(-1)^{|x||z|}[x,[y,z]] +(-1)^{|y||x|}[y,[z,x]] +(-1)^{|z||y|}[z,[x,y]]=0. \label{liesuperjacobi}
 \begin{array}[c]{c}
\text{Super Jacobi} \\
\text{identity}
\end{array}
\end{align}
In skew-symmetric superalgebras, the super hom-Jacobi identity can be presented equivalently in the form of super Leibniz rule for the maps $ad_x=[x,.]:\mathfrak{g} \rightarrow \mathfrak{g},$
\begin{equation}
[x,[y,z]]=[[x,y],z]+(-1)^{|x||y|}[y,[x,z]]. \label{liesuperjacobi:leibniz}
\end{equation}
However, for general linear maps $\alpha$, the Hom-Lie superalgebras are substantially different from Lie superalgebras, as all algebraic structure properties, morphisms, classifications and deformations become  dependent fundamentally on the joint simultaneous structure and properties of both operations, the linear map $\alpha$ and the bilinear product $[.,.]$ linked in an intricate way via the $\alpha$-twisted super-Jacoby identity \eqref{homliesuperjacobi}.
\end{remark}

As for all hom-superalgebras, an even homomorphism $ \phi : \mathfrak{g} \to \mathfrak{g}'$ of the hom-Lie superalgebras (or hom-Leibniz superalgebras) $(\mathfrak{g},[.,.],\alpha)$ and $(\mathfrak{g}',[.,.]',\beta)$ is said to be a homomorphism of hom-Lie superalgebras (or hom-Leibniz superalgebras), if $\phi [u,v]=[\phi (u), \phi (v))]'$ and $\phi \circ \alpha = \beta \circ \phi$.
The hom-Lie superalgebras (or hom-Leibniz superalgebras) $(\mathfrak{g},[.,.],\alpha)$ and $(\mathfrak{g}',[.,.]',\beta)$ are isomorphic, if there is a hom-Lie superalgebra (or hom-Leibniz superalgebras) homomorphism $ \phi : \mathfrak{g} \to \mathfrak{g}'$ such that $\phi$ be bijective \cite{Makhlouf2010:paradigmnonasshomalghomsuper}. Hom-subalgebras of hom-Lie superalgebra $(\mathfrak{g},[.,.],\alpha)$ are defined as  $\mathds{Z}_{2}$-graded linear subspaces $I =  (I\cap \mathfrak{g}_{0}) \oplus (I\cap \mathfrak{g}_{1}) \subseteq \mathfrak{g}$ closed under both $\alpha$ and $[.,.]$, that is  $\alpha(I) \subseteq I$ and $[I,I] \subseteq I$. Hom-subalgebra $I$ is called a hom-ideal of the hom-Lie superalgebra $\mathfrak{g}$, if $[I,\mathfrak{g}] \subseteq I$, and notation $I \triangleleft \mathfrak{g}$ is used in this case.
In super skew-symmetric hom-superalgebras, and in particular in hom-Lie superalgebras, by super skew-symmetry \eqref{superskewsymmetry},
$[I,\mathfrak{g}] \subseteq I$ is equivalent to $[\mathfrak{g},I] \subseteq I$, since
\begin{align*}
\forall  y = & y_0+y_1 \in I,\ x=x_0+x_1 \in \mathfrak{g} = \mathfrak{g}_{0} \oplus \mathfrak{g}_{1}, y_0\in I_{0}, y_1\in I_{1}, x_0\in \mathfrak{g}_{0}, x_1\in \mathfrak{g}_{1}:\\
[x,y] & = \sum_{j,k\in \mathds{Z}_{2}}[x_j, y_k] \stackrel{\eqref{superskewsymmetry}}{=}
\sum_{j,k\in \mathds{Z}_{2}} (-(-1)^{|j||k|})[y_k,x_j] \\
      & = -[y_0, x_0] - [y_0, x_1] + [y_1, x_1] - [y_1, x_0] \\
      & = [-y_0, x_0+x_1] + [y_1, x_1 - x_0] \\
      & \quad \in [I,\mathfrak{g}] +[I,\mathfrak{g}]\subseteq I+I = I,
      \quad \mbox{when}\quad  [I,\mathfrak{g}] \subseteq I,  \\
[y,x] & = \sum_{j,k\in \mathds{Z}_{2}}[y_k, x_j] \stackrel{\eqref{superskewsymmetry}}{=}
\sum_{j,k\in \mathds{Z}_{2}} (-(-1)^{|k||j|})[x_j, y_k] \\
& -[x_0,y_0] - [x_1,y_0] + [x_1,y_1] - [x_0,y_1] \\
& = [x_0+x_1,-y_0] + [x_1 - x_0,y_1] \\
& \quad \in [\mathfrak{g},I] +[\mathfrak{g},I]\subseteq I+I = I,
\quad \mbox{when} \quad  [\mathfrak{g},I] \subseteq I.
\end{align*}
Thus, in hom-Lie superalgebras, all right or left hom-ideals are two-sided hom-ideals.

Hom-Lie subalgebra $I$ of a hom-Lie superalgebra is called commutative if $[I,I] = 0$.
If $I$ is not abelian, then $[x,y]\neq 0$ for some non-zero elements $x,y \in I$.

\begin{definition}[\cite{liu2013hom}]
The center of a hom-Lie superalgebra $\mathfrak{g}$ is defined as
$$C(\mathfrak{g})=\{x \in \mathfrak{g}:  [x,\mathfrak{g}]=0\}.$$
The centralizer of a hom-ideal $I$ in  a hom-Lie superalgebra $\mathfrak{g}$ is defined as
$$C_{\mathfrak{g}}(I)=\{x \in \mathfrak{g}: [x,I]=0\}.$$
\end{definition}

In any hom-Lie superalgebra $(\mathfrak{g}=\mathfrak{g}_{0}\oplus \mathfrak{g}_{1},[.,.],\alpha)$, the center is the centraliser of hom-ideal $\mathfrak{g}$ in $(\mathfrak{g},[.,.],\alpha)$, that is $C(\mathfrak{g})=C_{\mathfrak{g}}(\mathfrak{g})$.
For any hom-ideal $I$ the centralizer $C_{\mathfrak{g}}(I)$ is a supersubspace
$C_{\mathfrak{g}}(I) =  (C_{\mathfrak{g}}(I)\cap \mathfrak{g}_{0}) \oplus (C_{\mathfrak{g}}(I)\cap \mathfrak{g}_{1}) $, since
\begin{align*}
\forall y \in I,\ x=x_0+x_1 \in \mathfrak{g} & = \mathfrak{g}_{0} \oplus \mathfrak{g}_{1}, x_0\in \mathfrak{g}_{0}, x_1\in \mathfrak{g}_{1}: \\
[x,y] & = [x_0+x_1,y]=[x_0,y]+[x_1,y]=0 \Leftrightarrow \\
[x_0,y] & =-[x_1,y]=[-x_1,y] \in \mathfrak{g}_0 \cap \mathfrak{g}_1 \cap I= \{0\} \Leftrightarrow \\
[x_0,y] & =[x_1,y]=0 \Leftrightarrow \\
x_j &\in C_{\mathfrak{g}}(I)\cap \mathfrak{g}_j, \ j\in \mathds{Z}_{2}.
\end{align*}
In general,
$[C_{\mathfrak{g}}(I),C(\mathfrak{g})(I)]\subseteq C_{\mathfrak{g}}(I)$ and
$\alpha(C_{\mathfrak{g}}(I))\subseteq C_{\mathfrak{g}}(I)$ are not assured,
since the equality $[[x_1,x_2],y]=0$ is not necessarily implied by $[x_1,y]=0$ and $[x_2,y] = 0$,
and $[x,y]=0$ does not necessarily imply $[\alpha(x),y]=0$
for $x_1, x_2, x\in \mathfrak{g}$ and $y\in I$.

\begin{lemma} Let $(\mathfrak{g},[.,.],\alpha)$ be a hom-Lie superalgebra.
If $(\mathfrak{g},[.,.],\alpha)$ is a multiplicative hom-Lie superalgebra
with surjective $\alpha$, that is $\alpha([.,.])= [\alpha(.),\alpha(.)]$ and $\alpha(\mathfrak{g})=\mathfrak{g}$, then
the center $C(\mathfrak{g})$ is a commutative hom-ideal in $(\mathfrak{g},[.,.],\alpha)$.
\end{lemma}

\begin{proof} The hom-supersubspace $C(\mathfrak{g})=(C(\mathfrak{g})\cap \mathfrak{g}_{0}) \oplus (C(\mathfrak{g})\cap \mathfrak{g}_{1})$ of the hom-Lie superalgebra $(\mathfrak{g},[.,.],\alpha)$ is closed under $[.,.]$ and $\alpha$. Indeed,
$\alpha(C(\mathfrak{g}))\subseteq C(\mathfrak{g})$,
since the preimage set $\alpha^{-1}(y) \neq \emptyset $ of any $y\in \mathfrak{g}$ is non-empty by surjectivity of $\alpha$, and
\begin{align*}
\forall \ x \in C(\mathfrak{g}), y\in \mathfrak{g}:
[\alpha(x),y] =  [\alpha(x),\alpha(\alpha^{-1}(y))]= \alpha([x,\alpha^{-1}(y)])=\alpha(\{0\})=\{0\}.
\end{align*}
Moreover, $[C(\mathfrak{g}),C(\mathfrak{g})] = [C(\mathfrak{g}),\mathfrak{g}]=\{0\} \subseteq C(\mathfrak{g})$ by definition of the center.
Hence, $C(\mathfrak{g})$ is commutative hom-ideal.
\end{proof}

\begin{lemma} \label{lem.hom.ideal}
Let $(\mathfrak{g},[.,.],\alpha)$ be a multiplicative hom-Lie superalgebra, $(\alpha([.,.])= [\alpha(.),\alpha(.)])$. If $I$ is a hom-ideal $I$ in $(\mathfrak{g},[.,.],\alpha)$
such that $\alpha$ is surjective on $I$, that is $\alpha(I)=I$, then
\begin{enumerate}[label=\upshape{(\roman*)},leftmargin=30pt]
\item
$C_{\mathfrak{g}}(I)$ is a hom-ideal in hom-Lie superalgebra $(\mathfrak{g},[.,.],\alpha)$.
\item
$C(I) = C_{I}(I)$ is a commutative hom-ideal in the hom-Lie superalgebra $(I,[.,.]_I,\alpha_I)$, where
$[.,.]_I$ and $\alpha_I$ are restrictions of $[.,.]$ and $\alpha$ to $I$.
\item
If $(\mathfrak{g},[.,.],\alpha)$ is a multiplicative hom-Lie superalgebra
with surjective $\alpha$, that is $\alpha([.,.])= [\alpha(.),\alpha(.)]$ and $\alpha(\mathfrak{g})=\mathfrak{g}$, then
the center $C(\mathfrak{g})$ is a commutative hom-ideal in $(\mathfrak{g},[.,.],\alpha)$.
\end{enumerate}
\end{lemma}
\begin{proof}
For any hom-ideal $I$, the hom-supersubspace $C_{\mathfrak{g}}(I)=(C_{\mathfrak{g}}(I))\cap \mathfrak{g}_{0}) \oplus (C_{\mathfrak{g}}(I))\cap \mathfrak{g}_{1})$ of the hom-Lie superalgebra $(\mathfrak{g},[.,.],\alpha)$ is closed under $[.,.]$ if $\alpha(I)=I$, since by super hom-Jacobi identity \eqref{homliesuperjacobi:leibniz}, definition of the centralizer,
and the condition $I=\alpha(I)$ of surjectivity of the restriction of $\alpha$ on $I$,
\begin{align*}
& \forall \ x \in I\cap H(\mathfrak{g}),\ y, z \in C_{\mathfrak{g}}(I)\cap H(\mathfrak{g}): \\
& [x,y]=0,\ [\alpha(y),[x,z]] = [\alpha(y),0]=0, \Rightarrow \\
& [\alpha(x),[y,z]]=[[x,y],\alpha(z)]+(-1)^{|x||y|}[\alpha(y),[x,z]] = 0,
\Rightarrow \\
& [I,[C_{\mathfrak{g}}(I),C_{\mathfrak{g}}(I)]]\stackrel{\alpha(I)=I}{=}
[\alpha(I),[C_{\mathfrak{g}}(I),C_{\mathfrak{g}}(I)]]=\{0\}
\Rightarrow \\
& [C_{\mathfrak{g}}(I),C_{\mathfrak{g}}(I)] \subseteq C_{\mathfrak{g}}(I).
\end{align*}
The hom-supersubspace
$C_{\mathfrak{g}}(I)=(C_{\mathfrak{g}}(I))\cap \mathfrak{g}_{0}) \oplus (C_{\mathfrak{g}}(I))\cap \mathfrak{g}_{1})$ is closed under $\alpha$, since
definition of the centraliser, surjectivity $\alpha(I)=I$
of $\alpha$ on $I$ and multiplicativity of $\alpha$ yield
\begin{align*}
& \forall \ x \in C_{\mathfrak{g}}(I)\cap H(\mathfrak{g}): \\
& [\alpha(C_{\mathfrak{g}}(I)),I] = [\alpha(C_{\mathfrak{g}}(I)),\alpha(I)]= \alpha([C_{\mathfrak{g}}(I),I])
\in \alpha(\{0\})= \{0\} \Rightarrow
\alpha(C_{\mathfrak{g}}(I))\in C_{\mathfrak{g}}(I).
\end{align*}
Thus, $C_{\mathfrak{g}}(I)$ is a hom-supersubalgebra in the hom-superalgebra $(\mathfrak{g},[.,.],\alpha)$.
Moreover,
\begin{align*}
& \forall \ x \in I\cap H(\mathfrak{g}),\ y\in \mathfrak{g}\cap H(\mathfrak{g}),
z \in C_{\mathfrak{g}}(I)\cap H(\mathfrak{g}): \\
& [x,y] \in I,\ [\alpha(y),[x,z]] = [\alpha(y),0]=0, \Rightarrow \\
& [\alpha(x),[y,z]]=[[x,y],\alpha(z)]+(-1)^{|x||y|}[\alpha(y),[x,z]] \in I,
\Rightarrow \\
& [I,[\mathfrak{g},C_{\mathfrak{g}}(I)]]\stackrel{\alpha(I)=I}{=}
[\alpha(I),[\mathfrak{g},C_{\mathfrak{g}}(I)]]\in I
\Rightarrow
[\mathfrak{g},C_{\mathfrak{g}}(I)] \subseteq C_{\mathfrak{g}}(I).
\end{align*}
Hence, $C_{\mathfrak{g}}(I)$ is a hom-ideal.
\end{proof}

\begin{lemma}
Let $\mathfrak{g}$ be a multiplicative hom-Lie superalgebras with surjective $\alpha$, that is $\alpha([.,.])= [\alpha(.),\alpha(.)]$ and $\alpha(\mathfrak{g})=\mathfrak{g}$, then the quotient $\mathfrak{g}/[\mathfrak{g},\mathfrak{g}]$ is commutative. Moreover, $[\mathfrak{g},\mathfrak{g}]$ is the smallest hom-ideal with this property: if $\mathfrak{g}/I$ is commutative for some hom-ideal $ I \subset \mathfrak{g}$, then $[\mathfrak{g},\mathfrak{g}] \subset I$.
\end{lemma}
\begin{proof}
	
It is obvious that $[[\mathfrak{g},\mathfrak{g}],[\mathfrak{g},\mathfrak{g}]] \subseteq [\mathfrak{g},\mathfrak{g}]$. Since $\mathfrak{g}$ is multiplicative hom-Lie superalgebra with surjective $\alpha$, we have $\alpha ([\mathfrak{g},\mathfrak{g}]) = [\alpha(\mathfrak{g}),\alpha(\mathfrak{g})] \subseteq [\mathfrak{g},\mathfrak{g}]$. Since $[\mathfrak{g},\mathfrak{g}] \subset \mathfrak{g}$, then $[[\mathfrak{g},\mathfrak{g}],\mathfrak{g}] \subseteq [\mathfrak{g},\mathfrak{g}]$. Hence $[\mathfrak{g},\mathfrak{g}]$ is is a hom-ideal.
Let $\bar{x} = x+ [\mathfrak{g},\mathfrak{g}] \in \mathfrak{g}/[\mathfrak{g},\mathfrak{g}]$, $\bar{y} = y+ [\mathfrak{g},\mathfrak{g}] \in \mathfrak{g}/[\mathfrak{g},\mathfrak{g}]$, where $x, y \in \mathfrak{g}$. So,
\begin{align*}
[\bar{x}, \bar{y}] =[x+ [\mathfrak{g},\mathfrak{g}], y+ [\mathfrak{g},\mathfrak{g}]]
=[x,y] +[\mathfrak{g},\mathfrak{g}]
=\overline{[x,y]} = \bar{0} \in \mathfrak{g}/[\mathfrak{g},\mathfrak{g}].
\end{align*}
Thus $\mathfrak{g}/[\mathfrak{g},\mathfrak{g}]$ is commutative.
If $\mathfrak{g}/I$ is commutative, then $\bar{x} = x+ I \in \mathfrak{g}/I$, $\bar{y} = y+ I \in \mathfrak{g}/I$ commute, for all $x, y \in \mathfrak{g},$
$$[\bar{x}, \bar{y}] = [x+I, y+I] = [x,y]+I =\bar{0} \in \mathfrak{g}/I,$$
$$\Rightarrow [x,y]+I =I \Rightarrow [x,y] \in I.$$
This is true for all $x,y \in \mathfrak{g}$, thus, $[\mathfrak{g},\mathfrak{g}] \subset I.$
\end{proof}

We are going to need the following definition throughout the rest of the article.

\begin{definition}[\cite{AmmarMakhloufSaadaoui2013:CohlgHomLiesupqdefWittSup,LarssonSilvestrov2005:GradedquasiLiealgebras}]
	A representation of the hom-Lie superalgebra $(\mathfrak{g},[.,.],\alpha)$ on a $\mathds{Z}_{2}$-graded linear space $V = V_{\bar{0}} \oplus V_{\bar{1}} $ with respect to $\beta \in gl(V)_{\bar{0}}$ is an even linear map $\rho : \mathfrak{g} \to gl(V)$, such that for all homogeneous $x,y \in \mathcal{H}(\mathfrak{g})$,
	\begin{eqnarray*}
	\rho (\alpha(x)) \circ \beta &=& \beta \circ \rho(x), \\
	\rho ([x,y]) \circ \beta &=& \rho (\alpha(x)) \circ \rho(y) -(-1)^{|x||y|} \rho (\alpha(y)) \circ \rho(x).
	\end{eqnarray*}
	
A representation $V$ of $\mathfrak{g}$ is called irreducible or simple, if it has no nontrivial subrepresentations. Otherwise $V$ is called reducible.
\end{definition}

For any linear transformation $T: X\mapsto X$ of a set $X$, and any nonnegative integer $s$, the $s$-times composition is $ T^s = T \circ \dots \circ T \quad ($s$-times), \quad T^0 =Id, \quad T^1 = T,$ and if $T$ is invertible with inverse map $T^{-1}\mathfrak{g}\rightarrow \mathfrak{g}$, then $T^{-s} = T^{-1} \circ \dots \circ T^{-1} \quad (s-times).$

Next, we recall the notion of $\alpha^s$-derivations.

\begin{definition}[\cite{AmmarMakhloufSaadaoui2013:CohlgHomLiesupqdefWittSup}]
	Let $(\mathfrak{g},[.,.]_\mathfrak{g},\alpha)$ be a hom-Lie superalgebra. For any nonnegative integer s, we call $D \in (End(\mathfrak{g}))_i$, where $i \in \mathds{Z}_2$, an $\alpha^s$-derivation of the multiplicative hom-Lie superalgebra $(\mathfrak{g},[.,.]_\mathfrak{g},\alpha)$, if for all homogeneous $x,y \in \mathcal{H}(\mathfrak{g})$,
	\begin{eqnarray}
D \circ \alpha &=& \alpha \circ D, \\
D([x,y]_{\mathfrak{g}}) &=& [D(x), \alpha^s (y)]_{\mathfrak{g}} + (-1)^{|D| |x|} [\alpha^s (x), D(y)]_{\mathfrak{g}}.
	\end{eqnarray}
\end{definition}

Denote by $Der_{\alpha^s}(\mathfrak{g}) = (Der_{\alpha^s}(\mathfrak{g}))_0 \oplus (Der_{\alpha^s}(\mathfrak{g}))_1$ the set of all $\alpha^s$-derivations of the hom-Lie superalgebra $(\mathfrak{g},[.,.],\alpha)$, and
$$
Der(\mathfrak{g}) = \bigoplus_{s \ge -1} Der_{\alpha^s} (\mathfrak{g}).$$
For any $D \in Der(\mathfrak{g})$ and $D' \in Der(\mathfrak{g})$, define their commutator $[D,D']$ as
\begin{equation} \label{brac.der.}
[D,D']_{\mathcal{D}} =D \circ D' - (-1)^{|D||D'|} D' \circ D.
\end{equation}

\begin{lemma}[\cite{AmmarMakhloufSaadaoui2013:CohlgHomLiesupqdefWittSup}] \label{lem.hom.lie.der.}
Let $(\mathfrak{g},[.,.]_\mathfrak{g},\alpha)$ be a multiplicative hom-Lie superalgebra and consider on $Der(\mathfrak{g})$ the endomorphism $\tilde{\alpha}$ defined by $\tilde{\alpha}(D) = \alpha \circ D$, then $(Der(\mathfrak{g}),[.,.]_{\mathcal{D}},\tilde{\alpha})$ is a hom-Lie superalgebra where $[.,.]_{\mathcal{D}}$ is given by \eqref{brac.der.}.
\end{lemma}

For any $x \in \mathfrak{g}$ satisfying $\alpha(x)=x$, the mapping
$ad_s(x): \mathfrak{g} \to \mathfrak{g}$ defined for all $y \in \mathfrak{g}$ by
$
ad_s(x)(y) = [x, \alpha^s(y)]_{\mathfrak{g}},
$
is a $\alpha^{s+1}$-derivation, called an inner $\alpha^{s+1}$-derivation \cite{AmmarMakhloufSaadaoui2013:CohlgHomLiesupqdefWittSup}, and the set
$Inn_{\alpha^{s+1}}(\mathfrak{g}) = \{ [x, \alpha^{s}(.)]_{\mathfrak{g}} \mid \ x \in \mathfrak{g}, \  \alpha(x)=x \}$ is a linear space in $Der_{\alpha^{s+1}} (\mathfrak{g})$.

\section{Complete hom-Lie superalgebras}
\label{sec:complhomliesuperalg}
In \cite{armakanrazavicomalg2020:completehomliesuper}, the authors introduced the notion of a complete hom-Lie superalgebra and in this section we state some results about it.

\begin{definition}[\cite{armakanrazavicomalg2020:completehomliesuper}] \label{Def:complhomLiesupalg}
	Hom-Lie superalgebra $\mathfrak{g}$ is called a complete hom-Lie superalgebra if $\mathfrak{g}$ satisfies the following two conditions:
\begin{align}
		C(\mathfrak{g})&= 0, \label{Cgeq0}\\
		Der_{\alpha^{s+1}}(\mathfrak{g}) &= ad_s(\mathfrak{g}). \label{Deralphaspl1geqadsg}	
\end{align}
\end{definition}

\begin{remark}
Let $\mathfrak{g}_0$ be a complete Lie algebra, then it is not necessary that $\mathfrak{g}$ be a complete hom-Lie superalgebra.

Let $(\mathfrak{g}_0, <.,.>, \alpha)$ be a semisimple hom-Lie algebra, $\mathfrak{g}_1$ be a
finite-dimensional linear space and $\tilde{\alpha}: \mathfrak{g} = \mathfrak{g}_0 \oplus \mathfrak{g}_1 \to \mathfrak{g} = \mathfrak{g}_0 \oplus \mathfrak{g}_1$ be an even endomorphism such that $\tilde{\alpha}|_{\mathfrak{g}_0} = \alpha$. Then by {\rm \cite{AmmarMakhloufJA2010:homliesuperaladmsuperalg}}, $(\mathfrak{g},[.,.],\alpha)$ is a hom-Lie superalgebra such that $[x,y]=0$ for all $x \in \mathfrak{g}_1$, $y \in \mathfrak{g}$ and $[x,y]=<x,y>$ for all $x, y \in \mathfrak{g}_0$ where $<.,.>$ is bracket operation of the hom-Lie algebra $\mathfrak{g}_0$. Since $C(\mathfrak{g}) \ne 0$, $\mathfrak{g}$ is not complete hom-Lie superalgebra but $\mathfrak{g}_0$ is complete, that is $C(\mathfrak{g}_0)= 0$ and $Der_{\alpha^{s+1}}(\mathfrak{g}_0) = ad_s(\mathfrak{g}_0)$.
\end{remark}

\begin{definition}
A hom-Lie superalgebra $(\mathfrak{g},[.,.],\alpha)$ is called solvable if $\mathfrak{g}^n = 0$ for some $n \in \mathds{N}$, where $\mathfrak{g}^n$, the members of the derived series of $\mathfrak{g}$, are defined inductively,
$$\mathfrak{g}^1 = \mathfrak{g}, \quad \mathfrak{g}^n = [\mathfrak{g}^{n-1}, \mathfrak{g}^{n-1}], \quad  n>1.$$
\end{definition}

Note that any commutative hom-Lie superalgebra is solvabe and for a multiplicative hom-Lie superalgebra $\mathfrak{g}$, we have $\alpha(\mathfrak{g}^n) \subseteq \mathfrak{g}^n$ for any $n$.

The hom-Lie superalgebra $(\mathfrak{g},[.,.],\alpha)$ is called semisimple if it does not contain any non-trivial solvable hom-ideal.

Let $\mathfrak{g}$ be a hom-Lie superalgebra and let $\Phi$ be a bilinear form on $\mathfrak{g}$. Recall that $\Phi$ is called invariant if $\Phi([x,y],z)=\Phi(x,[y,z])$ for all $x, y, z \in \mathfrak{g}$.
The invariant bilinear form associated to the adjoint representation of $\mathfrak{g}$ is called the Killing form on $\mathfrak{g}$.

Now, we check conditions under which the completeness of $\mathfrak{g}_0$ and $\mathfrak{g}$ are equivalent.

\begin{theorem}
Let $\mathfrak{g} = \mathfrak{g}_0 \oplus \mathfrak{g}_1$ be a multiplicative hom-Lie superalgebra with surjective $\alpha$ on $\mathfrak{g}$ and $\mathfrak{g}_0$. If $\mathfrak{g}$ has the non-degenerate Killing form, then $\mathfrak{g}_0$ is a complete hom-Lie algebra and $\mathfrak{g}$ is a complete hom-Lie superalgebra.
\end{theorem}

\begin{proof}
We know $\mathfrak{g}$ has non-degenerate Killing form, thus $Der_{\alpha^{s+1}}(\mathfrak{g}) = ad_s(\mathfrak{g})$ by \cite{ArmakanFarhSilv2021:ndKilformsHomLiesuper}. Since $\alpha$ is surjective, $C(\mathfrak{g})$ is commutative hom-ideal, and so $C(\mathfrak{g})$ is solvable. Hence $C(\mathfrak{g})=0$. Thus $\mathfrak{g}$ is complete.
	Let $\Phi$ be a non-degenerate Killing form of $\mathfrak{g}$. Then the restriction of $\Phi$ to $\mathfrak{g}_0$ is the non-degenerate Killing form of $\mathfrak{g}_0$. Hence $\mathfrak{g}_0$ is semisimple hom-Lie algebra and $Der_{\alpha^{s+1}}(\mathfrak{g}_0) = ad_s(\mathfrak{g}_0)$. Since $\alpha$ is surjective, $C(\mathfrak{g}_0)$ is commutative and solvable, So $C(\mathfrak{g}_0)=0$. Therefore $\mathfrak{g}_0$ is complete hom-Lie algebra.
\end{proof}

\begin{proposition}\label{prop.exist.comp.hi}
Let $\mathfrak{g}$ be a multiplicative hom-Lie superalgebra and I be a complete hom-ideal of $\mathfrak{g}$ with surjective $\alpha$ on both $\mathfrak{g}$ and I. There exists a hom-ideal J such that $\mathfrak{g} = I \oplus J$.
\end{proposition}

\begin{proof}
Let $J = C_{\mathfrak{g}}(I)$. Then $C_{\mathfrak{g}}(I)$ is a hom-ideal of $\mathfrak{g}$ by Lemma \ref{lem.hom.ideal}. Since $I$ is hom-ideal, $ad_s(x) \in Der_{\alpha^{s+1}}(I)$, for all $x \in \mathfrak{g}$. Since $I$ is complete, $Der_{\alpha^{s+1}}(I) = ad_s(I)$, so there exists a $\alpha^{s+1}$-derivation $D$ in $Der_{\alpha^{s+1}}(I)$ such that $ad_s(x) = D$. Hence there exists $r \in I$ such that
	\begin{equation*}
	D(t) = ad_s(x)(t) = [x, \alpha^s(t)] = [r, \alpha^s(t)],
	\end{equation*}
	for any $t \in I$. Then $[x-r, \alpha^s(t)]=0$ and $x-r \in C_{\mathfrak{g}}(I) = J$. Thus $x=r+l$, for some $l \in J$.
	On the other hand, since $I$ is complete, $I \cap J = I \cap C_{\mathfrak{g}}(I) = C(I) = 0 $. Therefore $\mathfrak{g} = I \oplus J$.
\end{proof}

\begin{definition}
	Let $\mathfrak{g}$ be a hom-Lie superalgebra and $h(\mathfrak{g}) = \mathfrak{g} \oplus Der(\mathfrak{g})$. The even bilinear map $[.,.]_h : h(\mathfrak{g}) \times h(\mathfrak{g}) \to h(\mathfrak{g})$ and a linear map $\alpha_h : h(\mathfrak{g}) \to h(\mathfrak{g})$ are defined in $h(\mathfrak{g})$ by
	\begin{equation*}
	[x+D, y+E]_h = [x,y]_g + D(y) - (-1)^{|x||E|} E(x) + [D,E]_{\mathcal{D}},
	\end{equation*}
	\begin{equation*}
	\alpha_h(x+D) = \alpha(x) + \alpha \circ D,
	\end{equation*}
	where $x,y \in \mathfrak{g}$, $D, E \in Der(\mathfrak{g})$ and $[.,.]_{\mathcal{D}}$ is bracket in $Der(\mathfrak{g})$ given by \eqref{brac.der.}. With the above notation, $h(\mathfrak{g})$ is a hom-Lie superalgebra. We call $h(\mathfrak{g})$ a holomorph hom-Lie superalgebra.
\end{definition}

We know that $(Der(\mathfrak{g}),[.,.]_{\mathcal{D}},\tilde{\alpha})$ is hom-Lie superalgebra by Lemma \ref*{lem.hom.lie.der.}. So we have the following results.

\begin{lemma} \label{lem.Cg}
Let $\mathfrak{g}$ be a multiplicative hom-Lie superalgebra and $(h(\mathfrak{g}),[.,.]_h,\alpha_h)$ be holomorph hom-Lie superalgebra.
	\begin{enumerate}[label=\upshape{(\roman*)},leftmargin=30pt]
		\item \label{i:lem.Cg:homsupalgequiv}  If $C(\mathfrak{g})=0$, then
		$C(Der(\mathfrak{g}))=\{ D\in Der(\mathfrak{g})| [D,Der(\mathfrak{g})]_{\mathcal{D}}=0 \}=0$.
		\item \label{ii:lem.Cg:homsupalgequiv} $\mathfrak{g}$ is hom-ideal of $h(\mathfrak{g})$ and $h(\mathfrak{g}) / \mathfrak{g} \simeq Der(\mathfrak{g})$.
		\item \label{iii:lem.Cg:homsupalgequiv} $\mathfrak{g} \cap C_{h(\mathfrak{g})}(\mathfrak{g}) = C(\mathfrak{g})$.
	\end{enumerate}
\end{lemma}

\begin{proof}
	
	Let $D\in (Der_{\alpha^s}(\mathfrak{g}))_i$, $i\in \mathds{Z}_2$ and $D\in C(Der(\mathfrak{g}))$. Then
	$[D,Der(\mathfrak{g})]_{\mathcal{D}}=0.$ So $[D,ad_s(x)]_{\mathcal{D}}=0$, hence $[D,ad_s(x)]_{\mathcal{D}}(y)=0$, for all $x, y\in \mathfrak{g}.$ Thus
	\begin{align*}
	& D(ad_s(x)(y)) - (-1)^{|D||x|} ad_s(x)(D(y))=0, \Rightarrow \\
	& D([x,\alpha^s(y)]) - (-1)^{|D||x|} [x,\alpha^s(D(y))]=0,
	\Rightarrow \\
	& D([x,\alpha^s(y)]) = (-1)^{|D||x|} [x,\alpha^s(D(y))],
	\Rightarrow \\
	& [D(x),\alpha^{2s}(y)] + (-1)^{|D||x|}[x, \alpha^s(D(y))] = (-1)^{|D||x|} [x,\alpha^s(D(y))],
	\Rightarrow \\
	& [D(x),\alpha^{2s}(y)]=0 \stackrel{C(\mathfrak{g})=0}{\Rightarrow} D(x)=0 \Rightarrow D=0.
	\end{align*}
	Therefore $C(Der(\mathfrak{g}))=0.$
	Next, $\mathfrak{g} \triangleleft \mathfrak{g}$, so $\mathfrak{g}$ is a hom-ideal of $h(\mathfrak{g})$ and $h(\mathfrak{g}) / \mathfrak{g} \simeq Der(\mathfrak{g}).$
	Now, let $x\in \mathfrak{g}.$ Then
	$$x \in C(\mathfrak{g}) \Longleftrightarrow [x,\mathfrak{g}]_\mathfrak{g}=0 \Longleftrightarrow [x,\mathfrak{g}]_h=0 \Longleftrightarrow x\in \mathfrak{g} \cap C_{h(\mathfrak{g})}(\mathfrak{g}).$$
	Hence $\mathfrak{g} \cap C_{h(\mathfrak{g})}(\mathfrak{g}) = C(\mathfrak{g})$.
\end{proof}

Now, we state some equivalence conditions
for a hom-Lie superalgebra to be complete, by using the notion of holomorph hom-Lie superalgebras.

\begin{definition}
	Let $\mathfrak{g}, \mathfrak{h}$ be two hom-Lie superalgebras. We call $\mathfrak{e}$ an extension of the hom-Lie superalgebra $\mathfrak{g}$
	by $\mathfrak{h}$, if there exists a short exact sequence
	$$ 0 \to \mathfrak{h} \to \mathfrak{e} \to \mathfrak{g} \to 0 $$
	of hom-Lie superalgebras and their morphisms.
	\begin{enumerate}[label=\upshape{(\roman*)},leftmargin=30pt]
		\item  	An extension
		$ 0 \to \mathfrak{h} \stackrel{i}{\to} \mathfrak{e} \stackrel{p}{\to} \mathfrak{g} \to 0  $
		is called trivial extension if there exists an hom-ideal $I \subset \mathfrak{e}$ such that $\mathfrak{e} = Ker(p) \oplus I.$
		\item An extension
		$ 0 \to \mathfrak{h} \stackrel{i}{\to} \mathfrak{e} \stackrel{p}{\to} \mathfrak{g} \to 0  $
		is called splitting extension if there exists an hom-supersubspace $S \subset \mathfrak{e}$ such that $\mathfrak{e} = Ker(p) \oplus S.$
	\end{enumerate}
\end{definition}

\begin{theorem} \label{thm:homsupalgequiv}
	For a multiplicative hom-Lie superalgebra $(\mathfrak{g},[.,.],\alpha)$ with surjective $\alpha$, the following conditions are equivalent:
	\begin{enumerate}[label=\upshape{(\roman*)},leftmargin=30pt]
		\item \label{i:thm:homsupalgequiv} $\mathfrak{g}$ is a complete hom-Lie superalgebra\textup{;}
		\item \label{ii:thm:homsupalgequiv} any splitting extension $\mathfrak{e}$ by $\mathfrak{g}$ is a trivial extension and $\mathfrak{e} = \mathfrak{g} \oplus C_{\mathfrak{e}}(\mathfrak{g});$
		\item \label{iii:thm:homsupalgequiv} $h(\mathfrak{g}) = \mathfrak{g} \oplus C_{h(\mathfrak{g})}(\mathfrak{g})$.
	\end{enumerate}
\end{theorem}

\begin{proof}
Let $\mathfrak{e}$ be a splitting extension by $\mathfrak{g}$ and assume \ref{i:thm:homsupalgequiv} holds. Hence $\mathfrak{g} \triangleleft \mathfrak{e}$ and $C_{\mathfrak{e}}(\mathfrak{g}) \triangleleft \mathfrak{e}$. By \ref{i:thm:homsupalgequiv}, $C(\mathfrak{g})=0$, so $\mathfrak{g} \cap C_{\mathfrak{e}}(\mathfrak{g}) = 0$. Since $\mathfrak{g} \triangleleft \mathfrak{e}$, $ad_s(e)(\mathfrak{g}) \subset \mathfrak{g}$, for any $e \in \mathfrak{e}$. Then the restriction $ad_s(e)|_{\mathfrak{g}}$ is a derivation of $\mathfrak{g}$. Since $\mathfrak{g}$ is comlpete, thus $ad_s(e)|_{\mathfrak{g}}$ is a $\alpha^{s+1}$-derivation of $\mathfrak{g}$. We set $\pi(e) = ad_s(e)|_{\mathfrak{g}}$, for all $e \in \mathfrak{e}$. Since $Der_{\alpha^{s+1}}(\mathfrak{g}) = ad_s(\mathfrak{g}) \simeq \mathfrak{g}$, the map $\pi$ is a homomorphism from $\mathfrak{e}$ onto $Der_{\alpha^{s+1}}(\mathfrak{g})$ and $Ker(\pi) = C_{\mathfrak{e}}(\mathfrak{g})$. Thus $\mathfrak{e} = \mathfrak{g} \oplus Ker(\pi)$. Therefore $\mathfrak{e} = \mathfrak{g} \oplus C_{\mathfrak{e}}(\mathfrak{g})$.
	Suppose \ref{ii:thm:homsupalgequiv} holds, then \ref{iii:thm:homsupalgequiv} is obvious by setting $\mathfrak{e} = h(\mathfrak{g})$.
	Next, suppose \ref{iii:thm:homsupalgequiv} holds. By Lemma \ref{lem.Cg}, $C(\mathfrak{g}) = \mathfrak{g} \cap C_{h(\mathfrak{g})}(\mathfrak{g})$. From \ref{iii:thm:homsupalgequiv}, $C_{h(\mathfrak{g})}(\mathfrak{g}) \simeq h(\mathfrak{g})/\mathfrak{g}$. By Lemma \ref{lem.Cg}, $h(\mathfrak{g})/\mathfrak{g} \simeq Der_{\alpha^{s+1}}(\mathfrak{g}) \simeq C_{h(\mathfrak{g})}(\mathfrak{g}) \simeq \mathfrak{g}$. Since $C(\mathfrak{g}) = 0$, then $\mathfrak{g} \simeq ad_s(\mathfrak{g})$. Thus, $Der_{\alpha^{s+1}}(\mathfrak{g}) = ad_s(\mathfrak{g})$. Therefore, $\mathfrak{g}$ is a complete hom-Lie superalgebra.
\end{proof}

In the following theorem, we check the condition under which
the completeness of $\mathfrak{g}$ and its ideals are equivalent.

\begin{theorem}
Let $(\mathfrak{g},[.,.],\alpha)$ be a multiplicative hom-Lie superalgebra and $\mathfrak{g} = I \oplus J$, where $I$ and $J$ are hom-ideals and $\alpha$ is surjective on $\mathfrak{g}$, $I$ and $J$. Then
	\begin{enumerate}[label=\upshape{(\roman*)},leftmargin=30pt]
		\item \label{i:homLiesupidials} $C(\mathfrak{g}) = C(I) \oplus C(J);$
		\item \label{ii:homLiesupidials} if $C(\mathfrak{g}) =0$, then
		\begin{equation*}
		ad_s(\mathfrak{g}) = ad_s(I) \oplus ad_s(J),
		\end{equation*}
		\begin{equation*}
		Der_{\alpha^{s+1}}(\mathfrak{g}) = Der_{\alpha^{s+1}}(I) \oplus Der_{\alpha^{s+1}}(J);
		\end{equation*}
		\item \label{iii:homLiesupidials} $\mathfrak{g}$ is complete if and only if $I$ and $J$ are complete.
	\end{enumerate}
\end{theorem}

\begin{proof}
	    \ref{i:homLiesupidials} By Lemma \ref{lem.hom.ideal}, $C(I)$ and $C(J)$ are hom-ideals of $\mathfrak{g}$. Furthermore, $I \cap J = 0$, and so $C(I) \cap C(J) = 0$.
		Let $a+b \in C(I) \oplus C(J)$, where $a \in C(I)$ and $b \in C(J)$. Thus $[a,I]=0$ and $[b,J]=0$. Let $m+n \in I \oplus J = \mathfrak{g}$, where $m \in I$ and $n \in J$. Then
		\begin{equation*}
		[a+b , m+n] = [a+b , m] + [a+b , n] = [a,m] + [b,m] + [a,n] + [b,n] = 0,
		\end{equation*}
		since $a,m \in I$, $b,n \in J$ and $[b,m], [a,n] \in I \cap J =0$. Therefore $a+b \in C(\mathfrak{g})$ and $C(I) \oplus C(J) \subseteq C(\mathfrak{g})$.
		Let $x = m+n \in C(\mathfrak{g})$, where $m \in I$ and $n \in J$. Then $[x,\mathfrak{g}] = [m+n, \mathfrak{g}] = [m+n, I+J]=0$. Since $x \in C(\mathfrak{g})$ and $[n,I] \subseteq [J,I]=0$, then
		\begin{equation*}
		[m,I] = [x-n , I] = [x,I] - [n,I] = 0.
		\end{equation*}
	Hence $m \in C(I)$. In the same way, $n \in C(J)$. Thus $C(\mathfrak{g}) \subseteq C(I) \oplus C(J)$.

        \ref{ii:homLiesupidials} \ For $D \in Der_{\alpha^{s+1}}(I)$, we define an extended linear transformation on $\mathfrak{g}$ by the equality $D(m+n)=D(m)$, for $m \in I$ and $n \in J$.  So $D \in Der_{\alpha^{s+1}}(\mathfrak{g})$, $Der_{\alpha^{s+1}}(I) \subseteq Der_{\alpha^{s+1}}(\mathfrak{g})$ and $Der_{\alpha^{s+1}}(J) \subseteq Der_{\alpha^{s+1}}(\mathfrak{g})$.
		Let $m \in I_i$, $n \in J$ and $D \in (Der_{\alpha^{s+1}}(\mathfrak{g}))_j$, where $i, j \in \mathds{Z}_2$. Since $I, J$ are hom-ideals,
		\begin{equation*}
		[D(m), n] = D([m,n]) = [D(m), \alpha^{s+1}(n)] + (-1)^{ij} [\alpha^{s+1}(m), D(n)] \in I \cap J.
		\end{equation*}
		The equality $I \cap J = 0$ yields $[D(m), \alpha^{s+1}(n)] = [\alpha^{s+1}(m), D(n)] = 0$. Let $D(m) = m' + n'$, where $m' \in I$ and $n' \in C(J)$. Then
		\begin{equation*}
		[D(m), \alpha^{s+1}(n)] = [m' + n', \alpha^{s+1}(n)] = [m', \alpha^{s+1}(n)] + [n' , \alpha^{s+1}(n)] = 0.
		\end{equation*}
		By \ref{i:homLiesupidials}, $n' = 0$. Hence $D(m) = m' \in I$. Thus $D(I) \subseteq I$.
In the same way, $D(J) \subseteq J$.
		Let $D \in Der_{\alpha^{s+1}}(\mathfrak{g})$ and $m + n \in I + J$, where $m \in I$ and $n \in J$. We define $\alpha^{s+1}$-derivations $E$ and $F$ by setting
		\begin{equation*}
		E(m+n)=D(m),
		\end{equation*}
		\begin{equation*}
		F(m+n)=D(n).
		\end{equation*}
Clearly, $E \in Der_{\alpha^{s+1}}(I)$ and $F \in Der_{\alpha^{s+1}}(J)$. Then $D = E + F \in Der_{\alpha^{s+1}}(I) + Der_{\alpha^{s+1}}(J)$. Therefore $Der_{\alpha^{s+1}}(\mathfrak{g}) = Der_{\alpha^{s+1}}(I) \oplus Der_{\alpha^{s+1}}(J)$ as a linear space, since $Der_{\alpha^{s+1}}(I) \cap Der_{\alpha^{s+1}}(J) = 0$.
		Now we prove that $Der_{\alpha^{s+1}}(I)$ and $Der_{\alpha^{s+1}}(J)$ are hom-ideals of hom-Lie superalgebra $Der_{\alpha^{s+1}}(\mathfrak{g})$. Let $E \in (Der_{\alpha^{s+1}}(I))_i$, $F \in (Der_{\alpha^{s+1}}(\mathfrak{g}))_j$ and $n \in J$. Using the commutator of $\alpha^{s+1}$-derivations, which is defined in \cite{AmmarMakhloufSaadaoui2013:CohlgHomLiesupqdefWittSup}, we have
		\begin{equation*}
		[F,E](n) = (F \circ E)(n) - (-1)^{ij} (E \circ F)(n) = 0.
		\end{equation*}
		Thus $Der_{\alpha^{s+1}}(I)$ is hom-ideal of $Der_{\alpha^{s+1}}(\mathfrak{g})$. Similarly $Der_{\alpha^{s+1}}(J)$ is hom-ideal of $Der_{\alpha^{s+1}}(\mathfrak{g})$.

        \ref{iii:homLiesupidials} \ Let $\mathfrak{g}$ be complete. Then $C(\mathfrak{g}) = 0$ and $C(I) = C(J) = 0$ by \ref{i:homLiesupidials}. Using $ad_s(\mathfrak{g}) = Der_{\alpha^{s+1}}(\mathfrak{g})$ and statements \ref{i:homLiesupidials} and \ref{ii:homLiesupidials}, we have
		\begin{equation*}
		ad_s(I) \oplus ad_s(J) = Der_{\alpha^{s+1}}(I) \oplus Der_{\alpha^{s+1}}(J),
		\end{equation*}
and $ad_s(I) \subseteq Der_{\alpha^{s+1}}(I)$ and $ ad_s(J) \subseteq Der_{\alpha^{s+1}}(J)$ yield $$ad_s(I) = Der_{\alpha^{s+1}}(I) \mbox{ and } ad_s(J) = Der_{\alpha^{s+1}}(J).$$
Therefore $I$ and $J$ are complete hom-Lie superalgebras.
		Conversly, let $I$ and $J$ are complete, then $C(\mathfrak{g}) = C(I) \oplus C(J) = 0$, by \ref{i:homLiesupidials}, and $Der_{\alpha^{s+1}}(\mathfrak{g}) = Der_{\alpha^{s+1}}(I) \oplus Der_{\alpha^{s+1}}(J) = ad_s(I) \oplus ad_s(J) = ad_s(\mathfrak{g})$, by \ref{ii:homLiesupidials}.
\end{proof}

\begin{definition} \label{def.simply.com}
	Let $\mathfrak{g}$ be a complete hom-Lie superalgebra. If any non-trivial hom-ideal of $\mathfrak{g}$ is not complete, then $\mathfrak{g}$ is called a simply complete hom-Lie superalgebra.
\end{definition}

A simple and complete hom-Lie superalgebra is a simply complete hom-Lie superalgebra.
Now, we want to state the relation between simply complete hom-Lie superalgebras and indecomposable complete hom-Lie superalgebras.

\begin{theorem} \label{theorem.simpl.compl.}
Let $(\mathfrak{g},[.,.],\alpha)$ be a complete multiplicative hom-Lie superalgebra with surjective $\alpha$ on $\mathfrak{g}$.
	\begin{enumerate}[label=\upshape{(\roman*)},leftmargin=30pt]
		\item  \label{i:Thm:complhomLie} $\mathfrak{g}$ can be decomposed into the direct sum of simply complete hom-ideals.
		\item  \label{ii:Thm:complhomLie} $\mathfrak{g}$ is simply complete if and only if it is indecomposable.
	\end{enumerate}
\end{theorem}

\begin{proof}
\ref{i:Thm:complhomLie}\ If $\mathfrak{g}$ is simply complete, then \ref{i:Thm:complhomLie} holds. If $\mathfrak{g}$ is not simply complete, then by Proposition \ref{prop.exist.comp.hi}, there exists a nonzero minimal complete hom-ideal $I$ of $\mathfrak{g}$ such that $\mathfrak{g} = I \oplus C_{\mathfrak{g}}(I)$. Since a hom-ideal of $C_{\mathfrak{g}}(I)$ is also a hom-ideal of $\mathfrak{g}$, by continuing this method for $C_{\mathfrak{g}}(I)$, we reach to the decomposition of $\mathfrak{g}$ into the simply complete hom-ideals.

\ref{ii:Thm:complhomLie}\ If $\mathfrak{g}$ is simply complete, then it is indecomposable by \ref{i:Thm:complhomLie}. Conversely, if $\mathfrak{g}$ is indecomposable, then it has no non-trivial hom-ideals. Hence, $\mathfrak{g}$ is simply complete by Definition \ref{def.simply.com}.
\end{proof}

\begin{definition}
	A subspace $I$ of a hom-Lie superalgebra $\mathfrak{g}$ is called a characteristic hom-ideal of $\mathfrak{g}$, if $D(I) \subset I$ for all $D \in Der(\mathfrak{g})$.
\end{definition}

\begin{lemma}
	Let $(\mathfrak{g},[.,.],\alpha)$ be a multiplicative hom-Lie superalgebra, $I$ be a characteristic hom-ideal of $\mathfrak{g}$ and $\alpha$ is surjective on $\mathfrak{g}$ and $I$. Then $I$ is hom-ideal of $\mathfrak{g}$.
\end{lemma}

\begin{proof}
	Let $x, y\in I$, since $\alpha$ is surjective on $I$ and $ad_s(\mathfrak{g})$ is a $\alpha^{s+1}$-derivation, then
	$$ [x,y] \stackrel{\alpha(I)=I}{=} [x, \alpha^s(t)] = ad_s(x)(t) \in I,$$
	where $t \in I$ and $\alpha^s(t)=y.$ Thus $[I,I] \subseteq I.$
	Next, $\alpha(I) \subseteq I$, since $\alpha$ is surjective on $I$.
	Let $y \in I$ and $a \in \mathfrak{g}$. Then
	$$ [a,y] \stackrel{\alpha(I)=I}{=} [a, \alpha^s(t)] =  ad_s(a)(t) \in I, $$
	where $t \in I$ and $\alpha^s(t)=y.$ Thus $[\mathfrak{g},I] \subseteq I.$ Therefore $I$ is a hom-ideal of $\mathfrak{g}$.
\end{proof}

\begin{theorem}
	Let $(\mathfrak{g},[.,.],\alpha)$ be a multiplicative hom-Lie superalgebra with surjective $\alpha$, $C(\mathfrak{g})=0$ and $ad_s(\mathfrak{g})$ be a characteristic hom-ideal of $Der(\mathfrak{g})$. Then $Der(\mathfrak{g})$ is complete. Furthermore, if $\mathfrak{g}$ is indecomposable and $[\mathfrak{g} , \mathfrak{g}] = \mathfrak{g}$, then $Der(\mathfrak{g})$ is simply complete.
\end{theorem}

\begin{proof}
The Hom-Lie superalgebra $\mathfrak{g}$ has trivial center, so $\mathfrak{g} \simeq ad_s(\mathfrak{g})$. Let $\mathfrak{p} = Der(\mathfrak{g})$, then $\mathfrak{g} \triangleleft \mathfrak{p}$. Let $\mathfrak{q}$ be a splitting extension by $\mathfrak{p}$, that is $\mathfrak{p} \triangleleft \mathfrak{q}$. Hence for all $q \in \mathfrak{q}$, we have $ad_s(q) \in Der_{\alpha^{s+1}}(\mathfrak{p})$. $\mathfrak{g}$ is a characteristic hom-ideal of $\mathfrak{p}$, so there exists $p \in \mathfrak{p}$ such that $ad_s(\mathfrak{p})|_{\mathfrak{g}} = ad_s(\mathfrak{q})|_{\mathfrak{g}}$. Then $ad_s(p-q)|_{\mathfrak{g}} =0 $ and $p-q \in C_{\mathfrak{q}}(\mathfrak{g})$. Hence we have $\mathfrak{q} = \mathfrak{p} + C_{\mathfrak{q}}(\mathfrak{g})$. On the other hand, $\mathfrak{p} \cap C_{\mathfrak{q}}(\mathfrak{g}) = C_{\mathfrak{p}}(\mathfrak{g}) = 0$ and $\mathfrak{p} \triangleleft \mathfrak{q}$, thus $\mathfrak{q} = \mathfrak{p} \oplus C_{\mathfrak{q}}(\mathfrak{g})$. Hence $C_{\mathfrak{q}}(\mathfrak{g}) \subseteq C_{\mathfrak{q}}(\mathfrak{p})$ and we have $\mathfrak{q} = \mathfrak{p} \oplus C_{\mathfrak{q}}(\mathfrak{p})$. Therefore by Theorem \ref{thm:homsupalgequiv}, $\mathfrak{p} = Der(\mathfrak{g})$ is a complete hom-Lie superalgebra.
	Now, assume that $Der(\mathfrak{g})$ is not simply complete. So there exists a simply complete hom-ideal $I$. By Proposition \ref{prop.exist.comp.hi}, there exists a hom-ideal $J$ such that $\mathfrak{p} = I \oplus J$.
	For any $x , y \in \mathfrak{g}$, there exists $x_1, y_1 \in I$ and $x_2, y_2 \in J$ such that $x= x_1 + x_2$ and $y = y_1 + y_2$. Thus $[x,y] = [x_1 + x_2 , y] = [x_1 , y] + [x_2 , y]$ such that $[x_1 , y] \in I \cap \mathfrak{g}$ and $[x_2 , y] \in J \cap \mathfrak{g}$. Hence $\mathfrak{g} = [\mathfrak{g} , \mathfrak{g}] = (I \cap \mathfrak{g}) \oplus (J \cap \mathfrak{g})$. $\mathfrak{g}$ is indecomposable, then $I \cap \mathfrak{g} = 0$ or $J \cap \mathfrak{g} = 0$. Hence $\mathfrak{g} \subset J$ and $I \subset C_{\mathfrak{p}}(\mathfrak{g}) = 0$. Therefore by Theorem \ref{theorem.simpl.compl.}, $\mathfrak{p} = Der(\mathfrak{g})$ is  simply complete.
\end{proof}


\section*{Acknowledgment}
The authors would like to thank Shiraz University, for financial support, which leads to the formation of this manuscript. This research is supported by grant no. 98GRC1M82582 Shiraz University, Shiraz, Iran.
\section*{Competing interests}
The authors declare that they have no competing interests.

\end{document}